.
\font\sets=msbm10.
\font\script=eusm10.
\font\stampatello=cmcsc10.
\font\stars=pzdr.

\def\spaziolungo{\qquad \qquad \qquad \qquad \qquad \qquad }
\def\1{{\bf 1}}

\def\avesum{\sum_{x\sim N}}

\def\square{\hbox{\vrule\vbox{\hrule\phantom{s}\hrule}\vrule}}
\def\defineq{\buildrel{def}\over{=}}

\def\doublesum{\mathop{\sum\sum}}

\def\integrale{\mathop{\int}}

\def\C{\hbox{\sets C}}
\def\N{\hbox{\sets N}}

\def\R{\hbox{\sets R}}
\def\Z{\hbox{\sets Z}}
\def\Corr{\hbox{\script C}}

\def\modSel{{\widetilde{J}}}
\def\Res{\mathop{{\rm Res}\,}}

\def\stella{\hbox{\stars I}}
\def\sinc{\hbox{\rm sinc}}

\par
\centerline{\bf A generalization of Gallagher's lemma for exponential sums}
\bigskip
\par
\centerline{\stampatello g. coppola\footnote{$^1$}{\rm titolare \enspace di \enspace un \enspace Assegno \enspace \lq \lq Ing. Giorgio Schirillo\rq \rq \thinspace \enspace \enspace dell'Istituto \enspace Nazionale \enspace di \enspace Alta \enspace Matematica \enspace (first author is supported by an INdAM Research Grant)} - m. laporta}
\bigskip
{

{\par
{\bf Abstract.} First we generalize a famous lemma of Gallagher 
on the mean square estimate for exponential sums
by plugging a weight
in the right hand side of Gallagher's original inequality. Then we apply it in the special case of the Ces\`aro weight
in order to establish some results mainly concerning the classical Dirichlet polynomials and the Selberg integrals of an arithmetic function $f$, that are tools for studying the distribution of $f$ in short intervals. Furthermore, we describe the smoothing process via self-convolutions of a weight that is involved into our Gallagher type inequalities, and compare it with 
the analogous process via the so-called correlations. Finally, we discuss a comparison argument in view of refinements on the Gallagher weighted inequalities according to different instances of the weight. 

\footnote{}{\par \noindent {\it Mathematics Subject Classification} (MSC 2010) : Primary $11{\rm L}07$; Secondary $11{\rm N}37$}
}
\bigskip
\par
\centerline{\bf 1. Introduction and statement of the main results.} 
\smallskip
\par
\noindent
In 1970 Gallagher ([G], Lemma 1) gave a general mean value estimate for series of the type
$$
S(t)\defineq\sum_{\nu} s(\nu)e(\nu t)\ ,
$$
\par
\noindent
where $e(x)\defineq e^{2\pi ix}$ as usual, 
the frequencies $\nu$ run over a (strictly increasing) sequence of real numbers 
and the coefficients $s(\nu)$ are complex numbers. 
More precisely, if $S(t)$ converges absolutely and 
$\delta, T>0$ are real numbers such that
$\theta\defineq \delta T\in(0,1)$, then 
$$
\|S\|^2_{2,T}\ll_{\theta}\delta^{-2}\int_{\R}\Big|\sum_{x<\nu \le x+\delta}s(\nu)\Big|^2dx\ ,
\leqno{(\star)}
$$
\par
\noindent
where for brevity we have written 
$$
\|S\|_{2,T}\defineq\|S\|_{L^2(-T,T)}=\left(\int_{-T}^T|S(t)|^2dt\right)^{1/2}\ .
$$
Hereafter, 
$A\ll_{\theta} B$ stands for $|A|\le CB$, where $C>0$ is an unspecified constant that depends on $\theta$. Typically, in the present context the bounds hold for $T\to 0$ or $T\to \infty$ (more precisely, for $|T|\le T_0$ with a sufficiently small $T_0>0$ or 
 for $T>T_0$ with a sufficiently large $T_0>1$). 
\par
We refer the reader to [H] for a first introduction to the large-sieve results, which constitute the main context 
where the inequality $(\star)$ has been widely applied (beware that in [H], Part III, $(\star)$ is referred to as {\it  Gallagher's second Lemma}). 
In particular, an immediate and renowned consequence is the so-called {\it Gallagher's Lemma for the Dirichlet series}
(see [G], Theorem 1).
Indeed, since an absolutely convergent Dirichlet series can be written as
$$
D(t)\defineq\sum_n a_nn^{it}=\sum_{\nu} s(\nu)e(\nu t)\ 
$$
\par
\noindent
by taking $\nu=(2\pi)^{-1}\log n$ and $s(\nu)=a_{e(\nu/i)}$, then, making  the substitution $x=\theta\log y$ in $(\star)$ with $\theta=(2\pi)^{-1}$ and recalling that $T=\theta\delta^{-1}$, one immediately has
$$
\|D\|^2_{2,T}\ll 
T^2 \int_0^{+\infty}\Big|\sum_{y<n\le ye^{1/T}}a_n\Big|^2{dy\over y}\ .
\leqno{(\stella)}
$$
\par
\noindent                
Inequalities of this type are important
tools for studying large values of the so-called Dirichlet polynomials, i.e.
finite Dirichlet series; 
in fact, they have large values frequencies strictly related to the behavior 
of their {\it moments}: see Ch.9 of [IK] and Ch.7 of [M] for basic knowledge on this topic. 
\bigskip
\par
\noindent
Here we give a more general version of $(\star)$ by plugging  a weight into the right hand side of it. To this end, for any {\it weight} $w:\R \rightarrow \C$ 
we set $w_\delta(x)=w(x)$ or $0$ according as $|x|\le \delta$ or not. The Fourier transform of $w_\delta$
is denoted by
$$
\widehat{w_\delta}(y)\defineq\int_{\R}w_\delta(t)e(-ty)dt=\int_{-\delta}^\delta w(t)e(-ty)dt\ .
$$
\bigskip
\par				
\noindent
Moreover, we recall that $L^1_{loc}(\R)$ denotes the space of all the locally summable functions on $\R$. 
\bigskip
\par				
\noindent
{\bf Lemma.} {\it Let the real numbers $\delta, T>0$ be fixed. Assume that, given
a (strictly increasing) sequence of real numbers $\nu$ 
and complex coefficients $s(\nu)$, the series \thinspace $S(t)=\sum_{\nu}s(\nu)e(\nu t)$ \thinspace is absolutely convergent. Then for every weight 
$w\in L^1_{loc}(\R)$ one has
$$
m_{\delta,T}\|S\|^2_{2,T}\le 
\int_{\R}\Big|\sum_{\nu}s(\nu)w_\delta(x-\nu)\Big|^2dx\ , 
\leqno{(\star\star)}
$$
\par
\noindent
where \enspace 
$\displaystyle{m_{\delta,T}\defineq \min_{|t|\leq T}|\widehat{w_\delta}(t)|^2}$. 
More in general, if $w\in L^1(\R)$, then
$$
m_{T}\|S\|^2_{2,T}\le 
\int_{\R}\Big|\sum_{\nu}s(\nu)w(x-\nu)\Big|^2dx\ , 
\leqno{(\star\!\star\!\star)}
$$
\par
\noindent
where \enspace 
$\displaystyle{m_{T}\defineq \min_{|t|\leq T}|\widehat{w}(t)|^2}$. 
}
\smallskip
\par
\noindent 
Though the proof of
$(\star\star)$ closely parallels Gallagher's original one, 
it is fully provided in \S2 together with the proof of
$(\star\!\star\!\star)$, while some further aspects of the Lemma are discussed in \S4. 
\smallskip
\par
\noindent
{\bf Remarks.}
\par
\noindent {\bf 1.} If $m_{\delta,T}=0$, as when $w_\delta$ is odd,
then $(\star\star)$ is trivial. Now, let $w_\delta$ be an even weight whose {\it self-convolution}
$$
(w_\delta\ast w_\delta)(x)\defineq \int_{\R}w_\delta(t)w_\delta(x-t)dt 
$$
is nonnegative. If $|s(\nu)|\le b(\nu)$, then $(\star\star)$ yields 
$$
m_{\delta,T}\int_{-T}^T\Big|\sum_\nu s(\nu)e(\nu y)\Big|^2dy\le 
\int_{\R}|\sum_{\nu} b(\nu)w_\delta(x-\nu)|^2dx\ .
$$
\par
\noindent 
Indeed, by applying the inequality $s(\nu_1)\overline{s(\nu_2)}+s(\nu_2)\overline{s(\nu_1)}\le 2|s(\nu_1)||s(\nu_2)|$ one has
$$
\int_{\R}\Big|\sum_{\nu}s(\nu)w_\delta(x-\nu)\Big|^2dx=
\sum_{\nu_1,\nu_2}s(\nu_1)\overline{s(\nu_2)}(w_\delta\ast w_\delta)(\nu_2-\nu_1)\le
$$
$$
\le
\sum_{\nu_1,\nu_2}b(\nu_1)b(\nu_2)(w_\delta\ast w_\delta)(\nu_2-\nu_1)=\int_{\R}|\sum_{\nu} b(\nu)w_\delta(x-\nu)|^2dx\ .
$$
\par
\noindent 
Note that, by an application of the Hardy-Littlewood majorant principle (see [M], Ch. 7, or  [L]) combined  with $(\star\star)$ for $b(\nu)$,
one would get the factor $3$ in the right hand side of
the above inequality .
\vfill\eject
 
\par                         
\noindent {\bf 2.}  Gallagher's original $(\star)$ is recovered from $(\star\star)$ by taking the weight $\delta^{-1}u_{\delta}$ associated to the {\it unit step} function
\smallskip
\par
\centerline{$
u(t)\defineq\cases{1\ \hbox{if}\ t>0,\cr 0\ \hbox{otherwise},\cr}$
with $\delta=\theta T^{-1}$ for $T>0$ and $\theta\in(0,1)$ fixed.} 
\smallskip
\noindent
In the following, we exploit the particular case of 
the {\it Ces\`aro weight}, 
$$
C_\delta(y)\defineq\max(1-\delta^{-1}|y|,0)\ ,
$$
\par
\noindent 
for which $(\star\star)$ yields the asymptotic inequality
$$
\|S\|^2_{2,T}\ll_\theta \delta^{-2}
\int_{\R}|\sum_{\nu}s(\nu)C_{\delta}(\nu-x)|^2dx
\ ,
\leqno{(\widetilde{\star})}
$$
\par
\noindent                                   
that was already established in the unpublished manuscript [CL1].
More explicitly, recalling that the well-known Fourier transform of $C_{\delta}$ is
(see also \S3)
$$
\widehat{C_{\delta}}(y)=
\cases{\displaystyle{
{{\sin^2(\pi\delta y)}\over {\pi^2 \delta y^2}}
}&  if $y\not=0$,\cr \cr
\delta&  if $y=0$,\cr}
$$
\par
\noindent 
the inequality $(\star\star)$ becomes
$$
\|S\|^2_{2,T}\le 
{{\pi^4 \theta^2 T^2}\over{\sin^4(\pi \theta)}}
\int_{\R}\Big|\sum_{|\nu-x|\le \theta/T}\Big(1-{T|\nu-x|\over\theta}\Big)s(\nu)\Big|^2dx\ .
\leqno{(\widetilde{\star\star})}
$$
\par
By applying $(\widetilde{\star})$, in \S2 we prove also the following modified version of $(\stella)$ for Dirichlet polynomials, but it may be easily generalized to absolutely convergent Dirichlet series. 
\smallskip
\par
\noindent {\bf Theorem 1}. {\it For every Dirichlet polynomial $\displaystyle{D(t)=\sum_n a_n n^{it}}$ one has, for} $T\to \infty$, 
$$
\|D\|^2_{2,T}\ll 
T^2 \int_{1}^{+\infty}\Big|\sum_{n}C_{y/T}(n-y)a_n\Big|^2{dy\over y}
+ \int_{1}^{+\infty}\Big(\sum_{y-\Delta\le n\le y+\Delta}|a_n|\Big)^2{dy\over y} \ ,
\leqno(\widetilde{\stella})
$$
\par
\noindent
{\it where} \enspace $\Delta=\Delta(y,T)\defineq y/T+O(y/T^2)$.
\medskip
\par
Now let us consider the special case of Dirichlet polynomials approximating Dirichlet series on the critical line $1/2+it$, namely
$$
P(t)\defineq \sum_{N_1\le n\le N_2}{{w(n)b(n)}\over {n^{1/2+it}}}\ , 
$$
where
\thinspace $N_1,N_2$ are positive integers,
$w$ is uniformly bounded and supported in \thinspace $[N_1,N_2]$, and
$b$ is an {\it essentially bounded} arithmetic function, that is
$|b(n)|\ll_{\varepsilon} n^{\varepsilon}$, $\forall \varepsilon>0$.
In \S2 we show the following consequence of Theorem 1 by applying $(\widetilde{\stella})$ with $a_n=w(n)b(n)n^{-1/2}$. 
\smallskip
\par
\noindent {\bf Corollary}. {\it For every $\varepsilon>0$, as $T\to \infty$, one has$$
\|P\|^2_{2,T}\ll_{\varepsilon} 
T^2 \int_{N_1/2}^{3N_2/2}\Big|\sum_{n}C_{y/T}(n-y){{w(n)b(n)}\over {n^{1/2}}}\Big|^2{dy\over y}
+ {{N_2^{1+\varepsilon}}\over {T^2}}\ . 
$$
}
\smallskip
Beyond the possible interest of our results within the general context 
of the exponential sums and
the possible further generalizations involving more Harmonic Analysis, here we 
wish to focus on another motivation for exploiting the particular case of  $(\star\star)$ with the Ces\`aro weight. 
First, let us recall that, taking inspiration from the classical method introduced by Selberg [S] to study 
the distribution of the prime numbers in short intervals, the so-called {\it Selberg integral} of an arithmetic function $f$ has been defined as (see [C1])
$$
\int_{h N^{\varepsilon}}^{N}\Big| \sum_{x<n\le x+h}f(n)-M_f(x,h)\Big|^2 dx\ , 
$$
\par
\noindent
where $N$ is an arbitrarily large integer, the real number $\varepsilon>0$ is arbitrarily small and
$M_f(x,h)$ is the expected mean value of the inner sum in the short interval $(x,x+h]$ (as $h=o(x)$ 
for $x\in[h N^{\varepsilon},N]$).
By a dyadic argument it is easily seen that 
the interval $[h N^{\varepsilon},N]$ can be replaced by $[N,2N]$ and that 
negligible remainder terms are generated when
the resulting integral on $[N,2N]$ is substituted by the discrete mean square (compare [CL]),
$$
J_f(N,h)\defineq \avesum \Big| \sum_{x<n\le x+h}f(n)-M_f(x,h)\Big|^2\, , 
$$
\par
\noindent
where $x\sim N$ means that $x\in(N,2N]\cap\N$. So that
one still refers to $J_f(N,h)$ as the {\it Selberg integral} of $f$.
\par				          
In this sense, if the coefficients $s(n)$ are assigned by an arithmetic function $s$ supported over $(N,2N]\cap\N$, then the right hand side of the Gallagher inequality $(\star)$ is de facto the Selberg integral of $s$, whenever we assume that 
$M_s(x,h)$ vanishes identically (in this case, we say that
$s$ is 
{\it balanced}). Indeed,
$$
J_s(N,\delta)=\avesum \Big| \sum_{x<n\le x+\delta}s(n)\Big|^2\ \hbox{with}\ \delta=o(N).
$$
\par
\noindent
Analogously, it transpires that the mean square
$$
\avesum \Big| \sum_{n}C_ \delta(n-x)s(n)\Big|^2
$$
\par
\noindent
emulates the integral on the right hand side of $(\widetilde{\star})$. 
Moreover, we showed that the inequality (see [CL], $\S4$)
$$
\avesum \Big| \sum_{n}C_ \delta(n-x)s(n)\Big|^2
\ll \avesum \Big| \sum_{x<n\le x+\delta}s(n)\Big|^2+\delta^3 \Vert s\Vert_{\infty}^2
$$
\par
\noindent
holds for every real and balanced function $s$ with  $\displaystyle{\Vert s\Vert_{\infty}\defineq \max_{N-\delta<n\le 2N+\delta}|s(n)|}$.
From this point of view, the inequality $(\widetilde{\star})$ can be proposed as a sort of refinement of Gallagher's inequality $(\star)$ (see \S3 for further discussions in this direction). 
\par
More in general, let us point out that the same mean value $M_f(x,h)$ appears in both $J_f(N,h)$ and the {\it modified} Selberg integral of $f$, i.e.
$$
\modSel_f(N,h)\defineq \avesum \Big| \sum_{n}C_h(n-x)f(n) - M_f(x,h)\Big|^2\ . 
$$
According to Ivi\'c [Iv], if the Dirichlet series $F(s)$ generated by
$f$ is meromorphic in $\C$ and absolutely convergent in the half-plane $\Re(s)>1$  at least, then the mean value 
takes the {\it analytic form}
$$
M_f(x,h)= hp_f(\log x),
$$
\par
\noindent
where $p_f(\log x)\defineq\Res_{s=1}F(s)x^{s-1}$ is the so-called {\it logarithmic polynomial} of $f$ (compare [CL]).

In [CL2] we exhibited the following {\it length-inertia} property for the Selberg integral, that allows to preserve non-trivial bounds as the length of the short interval 
increases: for $H\ge h$ with $h\to \infty$, $H=o(N)$ when $N\to \infty$ one has 
$$
J_f(N,H)\ll \left({H\over h}\right)^2 J_f(N,h) + J_f\left(N,H-h\left[{H\over h}\right]\right) + H^3 \left( \Vert f\Vert_{\infty}^2 + (\log N)^{2c}\right),
$$
\par				
\noindent
where $[\cdot]$ denotes the integer part,
$\displaystyle{\Vert f\Vert_{\infty}\defineq \max_{[N-H,2N+H]}|f|}$ and $c$ is the degree of the polynomial  $p_f$.
  
In order to establish
an analogous property for the modified Selberg integral of a real arithmetic function, we need to apply $(\widetilde{\star})$, while we underline that in [CL2] no Gallagher type inequality has been used in the proof of the  length-inertia  property for $J_f(N,H)$. Indeed, the last application of  $(\widetilde{\star})$ in the present paper is devoted to prove the following result (see \S2).
\smallskip
\par
\noindent {\bf Theorem 2.} {\it Let $f$ be a real arithmetic function for which the 
logarithmic polynomial $p_f(\log n)$ is defined. For $H\ge h$ with $h\to \infty$, $H=o(N)$ when $N\to \infty$ 
one has 
$$
\modSel_f(N,H)\ll H^2h^{-2} \modSel_f(N,h) 
              + \left(Nh^4H^{-2} 
               + H^3\right)  \Vert f\Vert_{\infty}^2 + H^3(\log N)^{2c}.
$$
}

\par                                       
After next section, that includes the proofs of the above results, in \S3 we first describe the smoothing process of a weight via self-convolutions and compare it with 
the analogous process performed trough the so-called correlations, which have been introduced in [CL]. Then we analyze the possible repercussions of such processes within the study of the weighted Selberg integral, with particular emphasis on the cases of the Ces\`aro weight and its relatives generated by iterations of the self-convolution. Finally, in \S4 a comparison argument is introduced in view of refinements on the right hand side of the Lemma inequalities, according to different instances of the weight $w$. 

\bigskip

\par
\centerline{\bf 2. Proofs of the main results.} 
\smallskip
\par
\noindent
{\stampatello Proof of the Lemma.}  First, note that there is nothing to prove
when the series 
$$
W_{\delta,s}(x)\defineq\sum_{\nu}s(\nu)w_\delta(x-\nu)
$$
is not square-integrable on $\R$, because the integral on the right hand side of 
$(\star\star)$ is not finite. Otherwise,
since by Lebesgue's dominated convergence theorem the Fourier transform of \thinspace $W_{\delta,s}$ \thinspace is 
$$
\widehat{W}_{\delta,s}(y)\defineq\int_{\R}W_{\delta,s}(x)e(-xy)dx=
\sum_{\nu} s(\nu)\int_{\R}w_\delta(x-\nu)e(-xy)dx=
$$
$$
=\sum_{\nu} s(\nu)e(-\nu y)\int_{\R}w_\delta(t)e(-ty)dt=S(-y)\widehat{w}_\delta(y)\ ,
$$
then $(\star\star)$ follows immediately by Plancherel's theorem
$$
\int_{\R}|W_{\delta,s}(x)|^2dx=\int_{\R}|\widehat{W}_{\delta,s}(y)|^2 dy
\ ,
$$
it being plain that
$$
\int_{\R}|S(-y)\widehat{w}_\delta(y)|^2 dy\ge m_{\delta,T}\int_{-T}^T|S(y)|^2 dy
\ .
$$
\smallskip

\noindent
Now, let us prove $(\star\!\star\!\star)$. We can clearly  assume that
$$
W_{s}(x)\defineq\sum_{\nu} s(\nu)w(x-\nu)\in L^2(\R)\ ,
$$
and apply $(\star\star)$  with  $\delta=n\in\N$ to write 
$$
m_{n,T}\int_{-T}^T|S(t)|^2dt\le 
\int_{\R}\Big|\sum_{\nu}s(\nu)w_n(x-\nu)\Big|^2dx\ .
$$
\par				
\noindent
Since Fourier transforms are continuous functions,  there exists 
$y_n\in[-T,T]$ such that 
$$
m_{n,T}=\min_{|y|\leq T}|\widehat{w}_n(y)|^2=|\widehat{w_n}(y_n)|^2=\Big|\int_{-n}^n w(t)e(-ty_n)dt\Big|^2\ .
$$
We can also assume that $y_n$ converges to some $y_0\in[-T,T]$ (this is surely true for some subsequence of $y_n$). 
Consequently, 
$w_n(t)e(-ty_n)$ converges to $w(t) e(-ty_0)$, while
it is plain that $|w_n(t)e(-ty_n)|\le |w(t)|$. Since
$w\in L^1(\R)$, then the dominated convergence theorem
yields 
$$
\lim_n m_{n,T}=\Big|\int_{\R} w(t)e(-ty_0)dt\Big|^2=|\widehat{w}(y_0)|^2\geq m_{T}\ .
$$                                                 
On the other side, the same theorem implies that
$$
\lim_n\int_{\R}\Big|\sum_{\nu}s(\nu)w_n(x-\nu)\Big|^2dx=
\int_{\R}\Big|\sum_{\nu}s(\nu)w(x-\nu)\Big|^2dx\ .
$$
Hence, $(\star\!\star\!\star)$ follows from the previous inequalities after passage to the limit as $n\to\infty$.\hfill $\square$

\bigskip

\par				
\noindent {\stampatello Proof of Theorem 1.} Let us apply $(\widetilde{\star})$ to
$$
D(t)=\sum_{\nu} s(\nu)e(\nu t),\enspace \hbox{with}\enspace \nu=(2\pi)^{-1}\log n\ ,\ s(\nu)=a_{e(\nu/i)}\ ,
$$
\par
\noindent
by taking $\theta=(2\pi)^{-1},\ T=\theta\delta^{-1}$ and $x=\theta\log y$. Thus, we get 
$$
\|D\|^2_{2,T}\ll \delta^{-2}\int_{\R}\Big|\sum_{|\nu-x|\le\delta}(1-|\nu-x|\delta^{-1})s(\nu)\Big|^2dx\ll 
$$
$$
\ll 
T^2 \int_0^{+\infty}\Big|\sum_{|\log n-\log y|\le 1/T}(1-T|\log n-\log y|)a_n\Big|^2{dy\over y}\ll 
$$
$$
\ll T^2 \int_{1}^{+\infty}\Big|\sum_{y-(1-1/\tau)y\le n\le y+(\tau-1)y}\Big(1-T\Big|\log\Big(1+{{n-y}\over y}\Big)\Big|\Big)a_n\Big|^2{dy\over y}\ , 
$$
\par
\noindent
where we have set $\tau\defineq e^{1/T}>1$ (note that $\tau\to 1$ as $T\to \infty$, so the $n-$sum is empty for $0<y<1$). 
\smallskip
\par
Since Taylor expansion yields 
$$
y-(1-1/\tau)y=y-{y\over T}+O\Big( {y\over {T^2}}\Big)\ , 
\quad 
y+(\tau-1)y=y+{y\over T}+O\Big( {y\over {T^2}}\Big)\ , 
$$
\par
\noindent
then the Ces\`aro weight, $\displaystyle{1-T|\log(1+(n-y)y^{-1})|}$, is bounded for the present range of $n$, while 
we have 
$$
1-T\Big|\log\Big(1+{{n-y}\over y}\Big)\Big|\ll {1\over T}
$$
\par
\noindent
in both ranges 
$$
0\le |n-(y-y/T)|\ll y/T^2
\enspace \hbox{\rm and} \enspace 
0\le |n-(y+y/T)|\ll y/T^2\ .
$$
\par
\noindent
Accordingly we write 
$$
\|D\|^2_{2,T}\ll 
T^2 \int_{1}^{+\infty}\Big|\sum_{y-y/T\le n\le y+y/T}\Big(1-T\Big|\log\Big(1+{{n-y}\over y}\Big)\Big|\Big)a_n\Big|^2{dy\over y} + 
$$
$$				
+ \int_{1}^{+\infty}\Big(\sum_{0\le |n-(y-y/T)|\ll y/T^2}|a_n|+\sum_{0\le |n-(y+y/T)|\ll y/T^2}|a_n|\Big)^2{dy\over y} \ .
$$
\par
\noindent
Then $(\widetilde{\stella})$ follows, since by Taylor expansion again we have 
$$
T\Big|\log\left(1+{{n-y}\over y}\right)\Big| - {{|n-y|}\over {y/T}}\ll 
{{T(n-y)^2}\over {y^2}}\ll 
{1\over T}
$$
\par
\noindent
for \enspace $y-y/T\le n\le y+y/T$ \enspace and this yields
$$
T^2\int_{1}^{+\infty}\Big|\sum_{y-y/T\le n\le y+y/T}\Big(1-T\Big|\log\Big(1+{{n-y}\over y}\Big)\Big|\Big)a_n\Big|^2{dy\over y}\ll 
$$
$$
\ll T^2 \int_{1}^{+\infty}\Big|\sum_{y-y/T\le n\le y+y/T}\Big(1-{{|n-y|}\over {y/T}}\Big)a_n\Big|^2{dy\over y} + 
 \int_{1}^{+\infty}\Big(\sum_{-y/T\le n-y\le y/T}|a_n|\Big)^2{dy\over y} \ ,
$$
\par				
\noindent
whence 
$$
\|D\|^2_{2,T}\ll T^2 \int_{1}^{+\infty}\Big|\sum_{y-y/T\le n\le y+y/T}\Big(1-{{|n-y|}\over {y/T}}\Big)a_n\Big|^2{dy\over y} + 
 \int_{1}^{+\infty}\Big(\sum_{-\Delta\le n-y\le \Delta}|a_n|\Big)^2{dy\over y} \ ,
$$
\par
\noindent
where recall that \enspace $\Delta=y/T+O(y/T^2)$\ .\hfill \square
\bigskip

\par
\noindent {\stampatello Proof of the Corollary.} 
Let us apply Theorem 1 to $D(-t)$ with $a_n=w(n)b(n)n^{-1/2}$
and write
$$
\|P\|^2_{2,T}\ll_{\varepsilon} 
T^2 \int_{1}^{+\infty}\Big|\sum_{y-y/T\le n\le y+y/T}\Big(1-{{|n-y|}\over {y/T}}\Big){{w(n)b(n)}\over {n^{1/2}}}\Big|^2{dy\over y}
+ N_2^{\varepsilon}\int_{1}^{+\infty}\Big(\sum_{y-\Delta\le n\le y+\Delta}{1\over {\sqrt n}}\Big)^2{dy\over y}\ . 
$$
\par
\noindent
Since $w$ has support  in $[N_1,N_2]$, then by using \enspace $\Delta=\Delta(y,T)\ll y/T$ \enspace one has
$$
\|P\|^2_{2,T}\ll_{\varepsilon} 
T^2 \int_{N_1/2}^{3N_2/2}\Big|\sum_{y-y/T\le n\le y+y/T}\Big(1-{{|n-y|}\over {y/T}}\Big){{w(n)b(n)}\over {n^{1/2}}}\Big|^2{dy\over y}
+ N_2^{\varepsilon}\int_{N_1/2}^{3N_2/2}{{\Delta(y,T)^2}\over {y^2}} dy\ll_{\varepsilon} 
$$
$$
\ll_{\varepsilon} T^2 \int_{N_1/2}^{3N_2/2}\Big|\sum_{y-y/T\le n\le y+y/T}\Big(1-{{|n-y|}\over {y/T}}\Big){{w(n)b(n)}\over {n^{1/2}}}\Big|^2{dy\over y}
+ N_2^{1+\varepsilon}{1\over {T^2}}\ .
$$
\par
\noindent
Hence the Corollary is proved. \hfill $\square$ 
\bigskip

\noindent
Before going to the proof of Theorem 2, let us prove an auxiliary proposition.
\smallskip
\par
\noindent {\bf Proposition 1.} {\it Let $w:\R \rightarrow \C$ be an uniformly bounded weight and let $f$ be an arithmetic function for which the 
logarithmic polynomial $p_f(\log n)$ is defined. For every fixed real number 
$\delta>0$ one has
$$
\avesum \Big|\sum_{n}w_\delta(n-x)f(n)-p_f(\log x)\sum_{n}w_\delta(n-x)\Big|^2
\ll \avesum\Big|\sum_{n}w_\delta(n-x)\tilde f(n)\Big|^2+
N^{-1}\delta^4(\log N)^{2c-2}\ ,
$$
where $\widetilde{f}(n)\defineq f(n)-p_f(\log n)$ assigns the {\it balanced} part of $f$, and $c$ is the degree of the polynomial $p_f$.
}
\par
\noindent {\stampatello Proof.} Since 
$$
\avesum \Big|\sum_{n}w_\delta(n-x)f(n)-p_f(\log x)\sum_{n}w_\delta(n-x)\Big|^2
\ll \avesum\Big|\sum_{n}w_\delta(n-x)\tilde f(n)\Big|^2+
$$
$$				
+\avesum \Big|\sum_{0\le |n-x|\le \delta}w(n-x)(p_f(\log n)-p_f(\log x))\Big|^2\ ,
$$
then the conclusion follows immediately from the mean value theorem applied to 
$p_f(\log n)-p_f(\log x)$. \hfill $\square$
\bigskip
\par
\noindent {\stampatello Proof of Theorem 2.} 
First, we recall here that for every $x\in\R$ one has
$$
p_f(\log x)\sum_{n}C_h(n-x)=hp_f(\log x)=M_f(x,h).
$$
For a balanced function $s$ supported over $(N,2N]\cap\N$, by taking
$\delta=h$ and $T=1/(2h)$ in 
$(\widetilde{\star\star})$ one sees that
$$
\int_{-{1\over {2h}}}^{1\over {2h}}\Big|\sum_{n\sim N}s(n)e(n\alpha)\Big|^2d\alpha\le 
{{\pi^4}\over{16h^2}}\int_{\R}|\sum_nC_h(n-x)s(n)|^2dx
\ll {1\over {h^2}}\modSel_s(N,h)+h\left\Vert s\right\Vert_{\infty}^2\ .
$$
\par
\noindent                                       
On the other side, by taking $w_\delta=C_\delta$ with $\delta=h\le H$ in the previous proposition and recalling that  $\tilde f$ is the balanced part of $f$, one gets
$$
\modSel_{f}(N,h)-\modSel_{\tilde f}(N,h)\ll N^{-1}h^4(\log N)^{2c-2}
\ll H^3 (\log N)^{2c}\ .
$$
\par
\noindent
Thus, without loss of generality we can assume that $f$ is balanced. Since by hypothesis $f$ is also real, we can apply the second formula\footnote{$^2$}{In [CL] formul\ae\ (59) are given 
in the context of a discussion about real and balanced functions that are supposed to be also essentially bounded.
However, it is easy to see that the latter assumption is in fact redundant for such formul\ae.} of (59) in [CL] and write
$$
\modSel_f(N,H)\ll {{1}\over {H^2}}
               \int_{-1/2}^{1/2}\left|{\cal F}_N(\alpha)\right|^2 \left|{\cal U}_H(\alpha)\right|^4d\alpha + H^3\left\Vert f\right\Vert_{\infty}^2\ ,
$$               
where
$$
{\cal U}_H(\alpha)\defineq\sum_{1\le n\le H}e(n\alpha)\ ,\qquad 
{\cal F}_N(\alpha)\defineq\sum_{n\sim N}f(n)e(n\alpha)\ .
$$
Let us recall that such exponential sums satisfy the well-known properties          
$$
{\cal U}_H(\alpha)\ll \min\Big(H,{1\over |\alpha|}\Big)\ \hbox{for}\ |\alpha|\le 1/2,\qquad 
\integrale_{|\alpha|\le 1/2}\left|{\cal F}_N(\alpha)\right|^2 d\alpha\ll
\sum_{n\sim N}f(n)^2\ .
$$
\par
\noindent
Now from all the previous inequalities it follows 
$$
 \int_{-1/2}^{1/2}\left|{\cal F}_N(\alpha)\right|^2 \left|{\cal U}_H(\alpha)\right|^4d\alpha
\ll
H^4 \int_{-1/(2h)}^{1/(2h)}\left|{\cal F}_N(\alpha)\right|^2 d\alpha 
        + h^4\integrale_{1/(2h)<|\alpha|\le 1/2}\left|{\cal F}_N(\alpha)\right|^2 d\alpha 
$$
$$
\ll
{{H^4}\over {h^2}}\left(\modSel_f(N,h)+h^3\left\Vert f\right\Vert_{\infty}^2\right)
        + h^4\integrale_{|\alpha|\le 1/2}\left|{\cal F}_N(\alpha)\right|^2 d\alpha 
\ll
{{H^4}\over {h^2}}\left(\modSel_f(N,h)+h^3\left\Vert f\right\Vert_{\infty}^2\right)+
h^4N\left\Vert f\right\Vert_{\infty}^2\ ,
$$
\par
\noindent
which implies the desired conclusion.\hfill \square

\vfill
\eject

\par				
\centerline{\bf 3. Two parallel ways of smoothing weights: self-convolution and autocorrelation.}
\smallskip
\par
\noindent
For every locally summable weight $w:\R \rightarrow \C$ and for some real number $\delta>0$  let us consider the {\it normalized}  self-convolution of $w_\delta$ given as
$$
\widetilde{w_{\delta}}(x)\defineq
{1\over {2\delta}}(w_{\delta} \ast w_{\delta})(x)={1\over {2\delta}}\int_{\R}w_\delta(t)w_\delta(x-t)dt\ .
$$
\par
\noindent
For example, the Ces\`aro weight $C_{\delta}$ is the normalized self-convolution of the restriction to $[-\delta/2,\delta/2]$ of $\1$, the
constantly 1 function: 
$$
C_{\delta}(x)
={1\over\delta}\!\!\!\integrale_{{|t|\le \delta/2}\atop{|x-t|\le \delta/2}}\!\!dt
={1\over {\delta}}(\1_{\delta/2}\ast \1_{\delta/2})(x)
=\!\widetilde{\enspace \1_{\delta/2}}(x)\ .
$$
\par
\noindent
It is well-known that the iteration of the self-convolution gives rise to
a process of smoothing (see [BN]). Moreover, the support of $\widetilde{w_{\delta}}$ is {\it doubled} with respect to the support of $w_\delta$
 in the sense that it is a subset of  $[-2\delta,2\delta]$. Because of 
the normalizing factor $(2\delta)^{-1}$, that takes into account the length of the integration interval, the magnitude of $w_\delta$ is not altered much
by the normalized self-convolution. More precisely, if one has $w_{\delta}\asymp 1$, i.e.
$1\ll w_{\delta}\ll 1$, in an interval of length $\gg \delta$, then there exists an interval of length $\gg \delta$ (not necessarily the same) where $\widetilde{w_{\delta}}\asymp 1$. From another well-known property of the convolution it follows that the Fourier transform of $\widetilde{w_{\delta}}$ is 
$$                                          
\widehat{\widetilde{w_{\delta}}}(y)=
{{\widehat{w_{\delta}}(y)^2}\over {2\delta}}\ .
$$ 
In particular, from
$$
\widehat{\1_{\delta}}(y)=\int_{-\delta}^{\delta}e(-ty)dt=2\delta\sinc(2\delta y)=
\cases{\displaystyle{
{{\sin(2\pi\delta y)}\over {\pi y}}}&  if $y\not=0$,\cr \cr
2\delta&  if $y=0$.\cr}
$$
we find that
$$
\widehat{C_{\delta}}(y)=\widehat{\widetilde{\enspace 1_{\delta/2}}}(y)=
{{\widehat{\1_{\delta/2}}(y)^2}\over {\delta}}=
\delta\sinc^2(\delta y)\ .
$$
\par
\noindent
Starting with 
$\1$, the normalized self-convolution generates recursively the family of 
Ces\`aro weights:
$$
C^{(j)}_{\delta}(x)\defineq \widetilde{C^{(j-1)}_{\delta/2}}(x), \ j\ge 1\ ,
$$
with the base steps
$
C^{(0)}_{\delta}(x)\defineq \1_{\delta}(x)$ and  $C^{(1)}_{\delta}(x)\defineq C_{\delta}(x)$. 
\par\noindent
Correspondingly, we have an inductive formula for the Fourier trasforms of these Ces\`aro weights.
\smallskip
\par
\noindent {\bf Proposition 2.} {\it For every $j\ge 0$ and every real number $\delta>0$,
$$
\widehat{C^{(j)}_\delta}(y)=
{{4\delta}\over{2^{2^j}}}\sinc^{2^j}\Big({{\delta y}\over{2^{j-1}}}\Big)\ .
$$
Consequently,   
$\widehat{C^{(j)}_\delta}(y)
\asymp_j \delta
$ at least for $0\le |y|\le 2^{j-2}\delta^{-1}$.
}
\smallskip

\noindent {\stampatello Proof.} The case $j=0$ is the above formula $\widehat{\1_{\delta}}$.  For $j\ge 1$, note that
$$
\widehat{C^{(j)}_\delta}(y)=\widehat{\widetilde{C^{(j-1)}_{\delta/2}}}(y)
=\delta^{-1}\widehat{C^{(j-1)}_{\delta/2}}(y)^2\ ,
$$
and more in general for $j\ge k\ge 0$, it is easily seen that
$$
\widehat{C^{(j)}_\delta}(y)={2^{a_k}\over \delta^{b_k}}\widehat{C^{(j-k)}_{\delta/2^k}}(y)^{2^k}\ ,
$$
\par
\noindent                     
where 
$$
a_k\defineq \sum_{i=1}^{k-1}i2^i=k2^k-2(2^k-1),\qquad b_k\defineq \sum_{i=0}^{k-1}2^i=2^k-1\ .
$$
By taking $k=j-1$ and setting $J=2^k=2^{j-1}$ for brevity, we obtain the stated formula
$$
\widehat{C^{(j)}_\delta}(y)=
{J^J\over (4\delta)^{J-1}}
\widehat{C_{\delta/J}}(y)^J={\delta\over{4^{J-1}}}
\sinc^{2J}\Big({\delta\over J}y\Big)\ .
$$
The remaining part of the statement is a straightforward consequence of this formula.\hfill $\square$
\bigskip

Such a process of 
{\it continuous} smoothing through the self-convolution of a weight $w_\delta$ has a {\it discrete} counterpart given by the {\it autocorrelation}\footnote{$^3$}{
Consistently with [CL], since no confusion can arise in the following, we will use the simpler term {\it correlation}.} of $w_\delta$:
$$
\Corr_{w_\delta}(a)\defineq \doublesum_{{n \thinspace \quad \thinspace m}\atop {n-m=a}}w_\delta(n)\overline{w_\delta(m)}\ . 
$$
\par
\noindent                                   
For example, since it turns out that
$$
C_\delta(t)
={1\over \delta} \sum_{a\le \delta-|t|}1 ={1\over \delta} \doublesum_{{a,b\le \delta}\atop {b-a=t}}1
={\Corr_{u_\delta}(t)\over \delta}\ ,
$$ 
\par
\noindent
then the Ces\`aro weight is the {\it normalized} correlation
of the unit step weight $u_\delta=\1_{[1,\delta]}$.
Note that the Ces\`aro weight is generated by both type of smoothing 
from a constantly one function. Moreover, through an iteration of the normalized correlation one might parallel the self-convolution process to generate the whole family of Ces\`aro weights $C^{(j)}_\delta$ with $j\ge 1$.
\par
An important aspect is that the exponential sums, whose coefficients are correlations of a weight $w$,  are non-negative. More precisely, 
$$
\sum_h \Corr_{w_\delta}(h)e(h\alpha)=\sum_h \doublesum_{n-m=h}w_\delta(n)\overline{w_\delta(m)}e(h\alpha) =
\Big|\sum_n w_\delta(n)e(n\alpha)\Big|^2\ . 
$$
\par
\noindent
A particularly well-known case is the {\it Fej\'er kernel} 
$$
\delta\sum_h C_\delta(h)e(h\alpha)=\sum_h \Corr_{u_\delta}(h)e(h\alpha)=
\Big|\sum_{1\le n\le \delta}e(n\alpha)\Big|^2\ .
$$
Such a {\it positivity} property is the complete analogous of the aforementioned fact that the Fourier 
transform of a self-convolution is a square. For this reason, in [CL] by abuse of notation we write 
$$
\widehat{w}_\delta(\alpha)=\sum_n w_\delta(n)e(n\alpha),\qquad
\widehat{\Corr_{w_\delta}}(\alpha)=\sum_h \Corr_{w_\delta}(h)e(h\alpha)\ , 
$$
\par
\noindent
and refer to such exponential sums as the {\it discrete Fourier transform} (DFT) of 
$w_\delta$ and $\Corr_{w_\delta}$, respectively.

By arguing as [CL] (compare formul\ae\ $(59)$) for the weighted Selberg integral $J_{w,f}$ of a real and balanced function $f$,
the positivity property provides the following alternative viewpoint of Gallagher's inequality:
$$
J_{w,f}(N,H)\defineq \avesum \Big|\sum_n w_H(n-x)f(n)\Big|^2
=\int_{-1/2}^{1/2}\Big|\sum_{n\sim N}f(n)e(n\alpha)\Big|^2 \left|\widehat{w}_H(\alpha)\right|^2d\alpha 
		+O\left(H^3\left\Vert f\right\Vert_{\infty}^2\right), 
$$
\par				
\noindent
where $\displaystyle{\Vert f\Vert_{\infty}\defineq \max_{N-H<n\le 2N+H}|f(n)|}$ and $\widehat{w}_H$ is the DFT of $w_H$. 
\smallskip
\par
\noindent
Trivially, this formula yields the inequality
$$
\min_{|y|\le 1/(2H)}\left|\widehat{w_H}(y)\right|^2\int_{-1/(2H)}^{1/(2H)}\Big|\sum_{n\sim N}f(n)e(n\alpha)\Big|^2 d\alpha 
\le J_{w,f}(N,H)+O\left(H^3\left\Vert f\right\Vert_{\infty}^2\right)\ . 
$$
\par
\noindent
In particular, recalling Proposition 2, we write
$$
\int_{-1/(2H)}^{1/(2H)}\Big|\sum_{n\sim N}f(n)e(n\alpha)\Big|^2 d\alpha 
\ll_j H^{-2}\tilde J_f^{(j)}(N,H)+H\left\Vert f\right\Vert_{\infty}^2\ , 
$$
\par
\noindent
where 
$$
J_f^{(j)}(N,H)\defineq \avesum \Big|\sum_n C^{(j)}_H(n-x)f(n)\Big|^2
$$ 
is what we might call the {\it $j$-th modified Selberg integral} of the balanced function $f$.
\smallskip

Such Gallagher type inequalities open the possibility of improving on the right hand side term by suitably picking a smoother weight from the Ces\`aro family, while 
Proposition 2 ensures that the order of magnitude of $\displaystyle{\min_{|y|\le 1/(2H)}\Big|\widehat{C^{(j)}_H}(y)\Big|^2}\asymp_jH^2$ is 
substantially unaffected by different choices of $j$. 
In this sense, our results of [CL] have already shed a light on the relation beween the cases $j=0$ (the Selberg integral) and $j=1$ (the modified Selberg integral). Finally, an additional feature offered by the discrete smoothing via correlations of a weight $w$ is a wider number theoretical perspective since
correlations can be plainly interpreted as a weighted count of solutions $n,m\sim N$ of the diophantine equation $n-m=a$.
\bigskip
\par				
\centerline{\bf 4. Comparing weights in view of Gallagher's generalized inequality.}
\smallskip
\par
\noindent
Bearing in mind the considerations of the previous section, here we compare weights in view of possible refinements of the right hand side of Gallagher's inequality
$(\star\star)$, here written as (see the proof of the Lemma)
$$
\int_{-T}^T|S(t)|^2dt\le m_{\delta,T}^{-1}
\int_{\R}|S(-y)\widehat{w}_\delta(y)|^2 dy
\ , 
$$
after assuming that $\displaystyle{m_{\delta,T}\defineq \min_{|t|\leq T}|\widehat{w_\delta}(t)|^2}\not=0$.
To this end, we give the following definition.
\smallskip
\par
\noindent
{\bf Definition.} Assume that $(\star\star)$ holds for both weights
$w_\delta, v_\delta$ such that 
\par
\centerline{$\displaystyle{m_{\delta,T}\defineq \min_{|t|\leq T}|\widehat{w_\delta}(t)|^2}\not=0$,\qquad
$\displaystyle{r_{\delta,T}\defineq \min_{|t|\leq T}|\widehat{v_\delta}(t)|^2}\not=0$.}
\par
\noindent
We say that $v_\delta$ is $T$-{\it better} than $w_\delta$ when
$$
{r_{\delta,T}\over m_{\delta,T}}
\ge {{|\widehat{v_\delta}(y)|^2}
\over {|\widehat{w_\delta}(y)|^2}}
 \quad \forall y\in\R\ .
$$
\smallskip
\par
\noindent
If so, it is plain that $v_\delta$ yields a refinement of
$(\star\star)$ with respect to $w_\delta$. Further, assuming that both weights are also positive, when $|y|\le T$ one has the following upper bound for the {\it gain}
$$
{
|\widehat{w_\delta}(y)|^2\over m_{\delta,T}}
-{{|\widehat{v_\delta}(y)|^2}
\over {r_{\delta,T}}}=
{{(|\widehat{w_\delta}(y)|^2-
m_{\delta,T})|\widehat{v_\delta}(y)|^2}\over
{m_{\delta,T}r_{\delta,T}}}
\leq 
\Big({||w_\delta||_1^2\over m_{\delta,T}}-
1\Big){||v_\delta||_1^2\over r_{\delta,T}}=
\Big({|\widehat{w_\delta}(0)|^2\over m_{\delta,T}}-
1\Big){|\widehat{v_\delta}(0)|^2\over r_{\delta,T}}
\ ,
$$  
\par				
\noindent
where $\|\cdot\|_1$ denotes the $L^1$ norm.
\smallskip
\par
\noindent
{\bf Example 1:} {\sl the Ces\`aro weights.}
Let us start by comparing $\1_\delta$ and  $C_\delta$.
Accordingly to the previous definition, $C_\delta$ is $T$-better than $\1_\delta$ if 
$$
{|\widehat{\1_\delta}(y)|^2\over|\widehat{C_\delta}(y)|^2}
\ge {\displaystyle{\min_{|t|\le T}|\widehat{\1_\delta}(t)|^2}\over \displaystyle{\min_{|t|\le T}|\widehat{C_\delta}(t)|^2}}\ \hbox{for every}\ y\in\R,\leqno(\1C)
$$
\par
\noindent
that is trivially true when $\widehat{\1_\delta}(y)=0$ for some $y\in[-T,T]$. Therefore, let us assume  
$\theta\defineq\delta T\in(0,1/2)$ for $\delta, T>0$, so that (see  $\S3$) 
$$
\min_{|t|\le T}|\widehat{\1_\delta}(t)|^2={\sin^2(2\pi\delta T)\over \pi^2 T^2}\not=0,\qquad
\min_{|t|\le T}|\widehat{C_\delta}(t)|^2=
{\sin^4(\pi\delta T)\over \pi^4T^4\delta^2}\not=0\ .
$$
\par
\noindent
Note that $(\1C)$ is true when $y=0$, since for $\delta T\in(0,1/2)$ one sees that
$$
{\displaystyle{\min_{|t|\le T}|\widehat{\1_\delta}(t)|^2}\over \displaystyle{\min_{|t|\le T}|\widehat{C_\delta}(t)|^2}}=
{4(\pi\delta T)^2\over \tan^2(\pi\delta T)}\le 4={|\widehat{\1_\delta}(0)|^2\over|\widehat{C_\delta}(0)|^2}\ .
$$
\par
\noindent
Thus, we consider the case $y\not=0$ and 
set $x\defineq y/T$ in $(\1C)$, that becomes 
$$
G_\theta(x)\defineq {|\widehat{C_\delta}(xT)|^2\over|\widehat{\1_\delta}(xT)|^2}=
{\tan^2(\pi\theta x)\over (2\pi\theta x)^2}
\le G_\theta(1)\ \hbox{for every}\ x\in\R\setminus\{0\}.
$$
It is easy to see that $G_\theta(x)$ satisfies the following properties:\par
(1) $G_\theta(x)$ is even with respect to both $x$ and $\theta$\par
(2) $G_\theta(x)$ is strictly increasing with respect to $x\in(0,1]$\par 
(3) $\displaystyle{\lim_{x\to 0} G_\theta(x)=1/4}$ $\forall\theta\in(0,1/2)$\par
(4) $G_\theta(1)$ is strictly increasing with respect to $\theta\in(0,1/2)$\par
(5) $\displaystyle{\lim_{\theta\to 0}G_\theta(1)=1/4}$, $\displaystyle{\lim_{\theta\to 1/2}G_\theta(1)=+\infty}$\par
(6) $G_\theta(x)=0\Longleftrightarrow x=k/\theta\ \forall k\in\Z\setminus\{0\}$\par
(note that $\forall k\in\Z\setminus\{0\}$ and $\forall\theta\in(0,1/2)$
one has $|k|/\theta> 2|k|\geq 2$)\par
(7) $G_\theta(x)\to+\infty$ as $x\to (2k+1)(2\theta)^{-1}\ \forall k\in\Z$\par
(note that $\forall k\in\Z$ and $\forall\theta\in(0,1/2)$
one has $|2k+1|(2\theta)^{-1}> |2k+1|\geq 1$)\par
\smallskip                                

\noindent
From properties (1)-(5) it follows that the above inequality for $G_\theta$
is true when $x\in[-1,1]\setminus\{0\}$, while properties (1), (6) and (7) imply that it is true for 
$|x|\in {\cal I}_{\theta,k}\defineq[(2k+1)(2\theta)^{-1}+\Delta'_k,(2k+3)(2\theta)^{-1}-\Delta''_k]\ \forall k\in\N$, where
both $\Delta'_k, \Delta''_k$ tend to $0^+$ as $k\to +\infty$. Since ${\cal I}_{\theta,k}$ tends to cover the whole interval $[2k+1,2k+3]$ as $\theta\rightarrow 1/2$ and $k\to +\infty$, then we deduce that $(\1C)$
is true for all $|y|\leq T$, and {\it almost everywhere} for $|y|>T$ as 
$\delta T\rightarrow 1/2$.
Thus, we say that
$C_\delta$ is {\it almost} $T$-better than $\1_\delta$ when 
$\delta T\rightarrow 1/2$. 

More in general, according to Proposition 2 of \S3, for every $j\ge 0$ one has
$$
{|\widehat{C^{(j+1)}_\delta}(y)|^2\over|\widehat{C^{(j)}_\delta}(y)|^2}= 
\left({\tan(\pi\delta y/2^j)\over\delta\pi y/2^{j-1}}\right)^{2^{j+1}}
=\left({\tan(\pi\theta x)\over2\pi\theta x}\right)^{2^{j+1}}=G_\theta(x)^{2^{j}}
\ ,
$$
where we have set 
$\theta=\delta T/2^j$, $x=y/T$ and $\displaystyle{
G_\theta(x)={\tan^2(\pi\theta x)\over (2\pi\theta x)^2}}$ as before.\par
Hence, we conclude that $C^{(j+1)}_\delta$ is {\it almost} $T$-better than $C^{(j)}_\delta$ for every $j\geq 0$, whenever
$\delta T\rightarrow 2^{j-1}$.
\smallskip

\noindent
{\bf Remark.} An effective use of $(\star\star)$ with $C^{(j)}_{\delta}$ requires finding explicit expressions of such weights. For example, the so-called  {\it Jackson-de La Vall\'e Poussin} weight $C^{(2)}_{\delta}$ (given by
the normalized self-convolution of  $C_{\delta/2}$) is the following cubic spline 
(see [BN], Problem 5.1.2 (v)):
$$
\delta^{-1}C_{\delta/2}\ast C_{\delta/2}(t)=\cases{\displaystyle{
{6|t|^3-6\delta t^2+\delta^3}\over{3\delta^3}}&  if $|t|\le \delta/2$,\cr \cr
\displaystyle{{2(\delta-|t|)^3\over 3\delta^3}}&  if $\delta/2<|t|\le \delta$,\cr\cr
0&  if $|t|>\delta$.\cr}
$$
\par
\noindent                  
Note that the support of  $C^{(2)}_{\delta}$ is $[-\delta,\delta]$, as expected. 
Evidently one could push forward the process by comparing arbitrary powers of $\widehat{\1_\delta}$. However,
the comparison between 
odd  and even powers seems to be cumbersome.
\medskip

\noindent
{\bf Example 2.} Given real numbers $\delta\ge\Delta>0$,
the {\it Lanczos} weight\footnote{$^4$}{In the literature, the Fourier transform   
 $\widehat{{\cal L}_{\delta,\Delta}}$ is known as {\it Lanczos kernel}.
 The diagram of ${\cal L}_{\delta,\Delta}$
 is an isosceles trapezium. Note that $\1_\delta\ge {\cal L}_{\delta,\Delta}\ge C_\delta={\cal L}_{\delta,\delta}$ for any $\delta\ge\Delta>0$.
 } is defined as (see [BN], Problem 5.1.2 (v))
$$
{\cal L}_{\delta,\Delta}(x)\defineq{1\over{\Delta}}(\1_{{\delta-\Delta/2}}\ast \1_{\Delta/2})(x)
=\cases{1&  if $|x|\leq \delta-\Delta$,\cr 
\displaystyle{{\delta-|x|}\over{\Delta}}&  if $\delta-\Delta<|x|\le \delta$,\cr
0&  if $|x|>\delta$,\cr
}
$$
whose Fourier transform is 
$$
\widehat{{\cal L}_{\delta,\Delta}}(y)
=(2\delta-\Delta)\sinc(\Delta y)\sinc((2\delta-\Delta) y)=
\cases{\displaystyle{
{{\sin(\pi\Delta y)\sin(2\pi\delta y-\pi\Delta y)}\over {(\pi y)^2\Delta}}}&  if $y\not=0$,\cr \cr
2\delta-\Delta&  if $y=0$.\cr}
$$
Since ${\cal L}_{\delta,\delta}=C_\delta$, then we can assume that 
$\delta>\Delta$. Let us
compare the weights ${\cal L}_{\delta,\Delta}$ and $\1_\delta$ by
verifying the inequality
$$
{|\widehat{\1_\delta}(y)|^2\over|\widehat{{\cal L}_{\delta,\Delta}}(y)|^2}
={{(2\delta)^2\sinc^2(2\delta y)}\over{
(2\delta-\Delta)^2\sinc^2(\Delta y)\sinc^2((2\delta-\Delta)y)}}
\ge {\displaystyle{\min_{|t|\le T}|\widehat{\1_\delta}(t)|^2}\over \displaystyle{\min_{|t|\le T}|\widehat{{\cal L}_{\delta,\Delta}}(t)|^2}}\ ,
$$
for every real $y$, and by assuming that $0<\Delta T< \delta T<1/2$. Then, it is easy to see that 
such inequality is satisfied by $y=0$. Therefore, for
$y\not=0$ we can write the left hand side as
$$
{|\widehat{\1_\delta}(y)|^2\over|\widehat{{\cal L}_{\delta,\Delta}}(y)|^2}={{(\Delta\pi y)^2\sin^2(2\delta\pi  y)}\over{
\sin^2(\Delta\pi y)\sin^2((2\delta-\Delta)\pi y)}}=
\Big({\Delta\pi y\over\tan((2\delta-\Delta)\pi y)}+{\Delta\pi y\over\tan(\Delta\pi y)}\Big)^2
\ .
$$
Since $0<\delta T<1/2$, then 
$$
\min_{|t|\le T}|\widehat{\1_\delta}(t)|^2=
{{\sin^2(2\pi\delta T)}\over{
(\pi T)^2}}\ ,
$$
which goes to $0$ as $\delta T\to 1/2$. Hence, we conclude that for $0<\Delta T< \delta T<1/2$ the Lanczos weight ${\cal L}_{\delta,\Delta}$ is {\it almost} $T$-better than $\1_\delta$ when 
$\delta T\rightarrow 1/2$. In a complete analogous way, we also see that, under the same conditions,
the Ces\`aro
weight $C_{\delta}$ is {\it almost} $T$-better than  ${\cal L}_{\delta,\Delta}$.\bigskip

\par
\centerline{\bf 5. Final considerations.} 
\smallskip
\par
\noindent
Because of the averaging over the inner short sum  coming from Ces\`aro weights, one could expect that, under suitable conditions on $f$, the modified Selberg integral $\modSel_f(N,h)$ should be more easily approachable than $J_f(N,h)$ (compare [CL], 
$\S0$). The process of smoothing described in \S3 and the comparison study in \S4 makes us to foretell that such a relaxing behavior might be  hopefully  shared by every $j$-th modified Selberg integral.  
Noteworthy, with the aid of $(\widetilde{\star})$ the first author [C3] has recently found a
way to deduce upper bounds for $J_f(N,h)$ from hypothetical estimates for $\modSel_f(N,h)$
by assuming that $f$ is balanced and essentially bounded. An application of such a method leads to a
non trivial estimate, $J_3(N,h)\ll N^{1+\varepsilon}h^{6/5}$,
for the Selberg integral
of the divisor function $d_3$, whenever one 
assumes that the sharp bound for the modified Selberg integral, $\modSel_3(N,h)
\ll N^{1+\varepsilon}h$, 
holds for every positive integer $h\ll N^{1/3}$ and for every real number $\varepsilon>0$.
This has to be compared with
the unconditional lower bound $Nh\log^4N\ll J_3(N,h)$ for $h\ll N^{1/3-\varepsilon}$ proved in [C2]. In [CL] we analyzed such conjecture on $\modSel_3(N,h)$
and analogous hypothesis on the modified Selberg integral $\modSel_k(N,h)$ of the divisor function
$$
d_k(n)\defineq \sum_{{n_1,\ldots,n_k}\atop {n_1\cdots n_k=n}}1\, ,\quad  (k\in \N)\ .
$$
\par				
\noindent
The importance of investigating about alternative ways of estimating the Selberg integral $J_k(N,h)$ of $d_k$ 
relies mainly on its strict connection with the 
$2k-$th moment of the Riemann zeta on the critical line (see [C1]). In particular, under the aforementioned conjectural bound of $\modSel_3(N,h)$, 
such an approach leads to the so-called {\it weak sixth moment for the Riemann zeta function} (see [CL], $\S8$). Finally, in the next future we are going to extend the results of the present paper to the general case of a weight $w$ that satisfies the hypothesis of the Lemma, mainly in view of possible applications to the study of the weighted Selberg integral $J_{w,f}(N,H)$ under suitable conditions on the function $f$.
\bigskip

\par
\noindent {\bf Acknowledgment}. The authors wish to thank Laura De Carli for having suggested us the reference [BN] and Alberto Perelli for helpful comments on an early version of [CL1].
\bigskip

\par
\centerline{\stampatello References}
\medskip
\item{\bf [BN]} P.L. Butzer and R.J. Nessel  \thinspace - \thinspace {\sl Fourier Analysis and Approximation} \thinspace - \thinspace Vol. I, Birkh\"auser Verlag, Basel und Stuttgart, 1971. 
\item{\bf [C1]} Coppola, G. \thinspace - \thinspace {\sl On the Selberg integral of the $k$-divisor function and the $2k$-th moment of the Riemann zeta-function} \thinspace - \thinspace Publ. Inst. Math. (Beograd) (N.S.) {\bf 88(102)} (2010), 99--110. $\underline{\tt MR\enspace 2011m\!:\!11173}$.  
\item{\bf [C2]} Coppola, G. \thinspace - \thinspace {\sl On some lower bounds of some symmetry integrals} \thinspace - \thinspace Afr. Mat. {\bf 25.1} (2014), 183--195. $\underline{\tt MR\enspace 3165958}$ (a draft version available at http://arxiv.org/abs/1003.4553, version 2).
\item{\bf [C3]} Coppola, G. \thinspace - \thinspace {\sl On the Selberg integral of the three-divisor function $d_3$} \thinspace - \thinspace available online at the address http://arxiv.org/abs/1207.0902, version 3.
\item{\bf [CL]} Coppola, G. and Laporta, M. \thinspace - \thinspace {\sl Generations of correlation averages} \thinspace - \thinspace Journal of Numbers, Volume 2014 (2014), Article ID 140840, 13 pages (a draft version available at http://arxiv.org/abs/1205.1706v3).
\item{\bf [CL1]} Coppola, G. and Laporta, M. \thinspace - \thinspace {\sl A modified Gallagher's Lemma} \thinspace - \thinspace unpublished \thinspace - \thinspace
(a draft version available at http://arxiv.org/abs/1301.0008v1).
\item{\bf [CL2]} Coppola, G. and Laporta, M. \thinspace - \thinspace {\sl Symmetry and short interval mean-squares} \thinspace - \thinspace submitted \thinspace - \thinspace (a draft version available at http://arxiv.org/abs/1312.5701).
\item{\bf [G]} Gallagher, P. X. \thinspace - \thinspace {\sl A large sieve density estimate near $\sigma =1$} \thinspace - \thinspace Invent. Math. {\bf 11} (1970), 329--339. $\underline{\tt MR\enspace 43\# 4775}$.
\item{\bf [H]} Huxley, M.N. \thinspace - \thinspace {\sl The Distribution of Prime Numbers} \thinspace - \thinspace Oxford University Press, Oxford (1972).
\item{\bf [Iv]} Ivi\'c, A. \thinspace - \thinspace {\sl On the mean square of the divisor function in short intervals} \thinspace - \thinspace J. Th\'eor. Nombres Bordeaux {\bf 21} (2009), no. {\bf 2}, 251--261. $\underline{\tt MR\thinspace 2010k\!:\!11151}$.
\item{\bf [IK]} Iwaniec H. and Kowalski E.\thinspace - \thinspace {\sl Analytic Number Theory} \thinspace - \thinspace AMS Colloquium Publications, {\bf 53}, Providence, RI, (2004). $\underline{\tt MR\thinspace 2005h\!:\!11005}$.
\item{\bf [L]} Laporta, M. \thinspace - \thinspace {\sl Some remarks on the Hardy-Littlewood majorant principle for exponential sums} \thinspace - \thinspace Ricerche Mat., (published online since August 1, 2014). 
\item{\bf [M]} Montgomery, H.L. \thinspace - \thinspace {\sl Ten Lectures on the Interface Between Analytic Number Theory 
and Harmonic Analysis} \thinspace - \thinspace CBMS Regional Conf. Ser. in Math. {\bf 84}, Amer. Math. Soc., Providence, RI, (1994). 
\item{\bf [S]} Selberg, A. \thinspace - \thinspace {\sl On the normal density of primes in small intervals, and the difference between consecutive primes} \thinspace - \thinspace Arch. Math. Naturvid. {\bf 47} (1943), {\bf no.6}, 87--105. $\underline{\tt MR\enspace 7,48e}$.

\bigskip

\par
\leftline{\tt Giovanni Coppola\spaziolungo \spaziolungo \qquad \qquad \enspace \thinspace Maurizio Laporta}
\leftline{\tt Universit\`a degli Studi di Salerno\spaziolungo \thinspace Universit\`a degli Studi di Napoli}
\leftline{\tt Home address \negthinspace : \negthinspace Via Partenio \negthinspace 12 \negthinspace -\spaziolungo Dipartimento di Matematica e Appl.}
\leftline{\tt - 83100, Avellino(AV), ITALY\spaziolungo \qquad \qquad \qquad \qquad \enspace \thinspace Compl.Monte S.Angelo}
\leftline{\tt e-page : $\! \! \! \! \! \!$ www.giovannicoppola.name\qquad \qquad \qquad \qquad \qquad \quad \thinspace Via Cinthia - 80126, Napoli, ITALY}
\leftline{\tt e-mail : gcoppola@diima.unisa.it\spaziolungo \qquad \enspace \thinspace e-mail : mlaporta@unina.it}

\bye